\documentclass{amsart}
\usepackage{amsmath,amsfonts,graphicx}
\usepackage[pdftex]{hyperref}

\pdfoutput=1

\begin{document}

\newtheorem{thm}{Theorem}[section]
\newtheorem{lem}[thm]{Lemma}
\newtheorem{cor}[thm]{Corollary}
\newtheorem{prop}[thm]{Proposition}
\newtheorem{conj}[thm]{Conjecture}
\newtheorem{mainlem}[thm]{Main Lemma}

\theoremstyle{definition}
\newtheorem{defn}{Definition}[section]

\theoremstyle{remark}
\newtheorem{rmk}{Remark}[section]

\def\square{\hfill${\vcenter{\vbox{\hrule height.4pt \hbox{\vrule width.4pt
height7pt \kern7pt \vrule width.4pt} \hrule height.4pt}}}$}

\newenvironment{pf}{{\it Proof:}\quad}{\square \vskip 12pt}
\newcommand{\dt}{\ensuremath{\text{det}}}
\newcommand{\h}{{\ensuremath{\text{hyp}}}}

\title{Minimal Volume Alexandrov Spaces}
\author{Peter A. Storm}\thanks{Supported by the NSF through VIGRE}

\date{September 18, 2002}

\begin{abstract}
Closed hyperbolic manifolds are proven to minimize volume over all Alexandrov spaces with curvature bounded below by $-1$ in the same bilipschitz class.  As a corollary compact convex cores with totally geodesic boundary are proven to minimize volume over all hyperbolic manifolds in the same bilipschitz class.  Also, closed hyperbolic manifolds minimize volume over all hyperbolic cone-manifolds in the same bilipschitz class with cone angles $\le 2\pi$.  The proof uses techniques developed by Besson-Courtois-Gallot.  In $3$ dimensions, this result provides a partial solution to a conjecture in Kleinian groups concerning acylindrical manifolds.
\end{abstract}
\maketitle

\section{Introduction}
\label{intro}
To state this paper's main result, let $N$ be a compact irreducible acylindrical $3$-manifold which admits a convex cocompact hyperbolic metric on its interior.  Then by Thurston's Geometrization and Mostow Rigidity there exists a convex cocompact hyperbolic manifold $M_0$ homeomorphic to int$(N)$ such that the convex core $C_{M_0}$ has totally geodesic boundary \cite[pg.14]{Th2}.  Here we prove

\vskip 6pt
\noindent\textbf{Theorem \ref{3-cor}.}  \itshape With the above notation, let $M$ be any hyperbolic manifold homotopy equivalent to $N$.  Then 
$$\text{Vol}(C_M) \ge \text{Vol}(C_{M_0}).$$
\normalfont
\vskip 6pt

This theorem is related to a conjecture in Kleinian groups made by Bonahon (see also \cite{CMT}).  To state things precisely, let $\mathcal{I}(N)$ denote the set of isometry classes of hyperbolic $3$-manifolds homotopy equivalent to $N$, and define a volume function $\text{Vol}: M \in \mathcal{I}(N) \longmapsto \text{Vol}(C_M).$  The conjecture states $M_0$ is the unique global minimum of Vol over $\mathcal{I} (N)$.  This paper proves $M_0$ is a global minimum of Vol.  Previous progress on this conjecture was made by Bonahon \cite{Bon}.  Using different techniques, and in a slightly more general setting, Bonahon proved $M_0$ is a strict local minimum of Vol.

The above theorem is proven by extending a minimal volume result of Besson-Courtois-Gallot \cite{BCGlong} to Alexandrov spaces with curvature bounded below by $-1$.
\vskip 6pt
\noindent\textbf{Theorem \ref{technical theorem}.}  \itshape Let $X$ be an Alexandrov space with curvature bounded below by $-1$, and $M_\h$ a closed  hyperbolic manifold.  If $X$ and $M_\h$ are bilipschitz, then
$$\text{Vol}(X) \ge \text{Vol}(M_\h).$$
\normalfont

\vskip 6pt
Theorem \ref{technical theorem} is used to study the invariant
$$\mathcal{V} (N) := \inf_{M \in \text{cc}_{\text{top}}(N)} \{ \text{volume of the convex core } C_M \text{ of } M \},$$
where $\text{cc}_{\text{top}} (N)$ is the set of isometry classes of complete convex cocompact hyperbolic manifolds diffeomorphic to the \textit{interior} of a smooth compact $n$-manifold $N$.  
Using the main theorem we prove

\vskip 6pt
\noindent\textbf{Theorem \ref{main2}.}  \itshape Let $N$ be a smooth  $n$-manifold.  If there exists an $M_0 \in \text{cc}_{\text{top}}(N)$ such that $\partial C_{M_0} \subset M_0$ is a totally geodesic submanifold, then 
$$\mathcal{V}(N) = \text{Vol} (C_{M_0}).$$
\normalfont

\vskip 6pt
{\noindent}As an immediate corollary, the Gromov norm of the doubled manifold $DN$ is related to $\mathcal{V}$ by the formula $\| [DN] \| v_n = 2 \mathcal{V} (N)$  (where $v_n$ is the volume of a regular ideal simplex in $\mathbb{H}^n$).  Applying Theorem \ref{main2} to $3$-manifolds yields Theorem \ref{3-cor}, which is stated above.

Theorem \ref{technical theorem} also proves two corollaries concerning cone-manifolds.
\vskip 6pt
\noindent\textbf{Theorem \ref{cone1}.}  \itshape Let $X$ be an $n$-dimensional cone-manifold with all cone angles $\le 2\pi$ and sectional curvatures $K \ge -1$ on smooth points.  Let $M_\h$ be a closed  hyperbolic $n$-manifold.  If $X$ and $M_\h$ are bilipschitz then
$\text{Vol} (X) \ge \text{Vol} (M_\h).$
\normalfont
\vskip 6pt

{\noindent}In the case when $n=3$, applying the Manifold Hauptvermutung yields the stronger statement
 
\vskip 6pt
\noindent\textbf{Corollary \ref{cone2}.}  \itshape 
Let $X$ be a $3$-dimensional cone-manifold with all cone angles $\le 2\pi$ and sectional curvatures $K \ge -1$ on smooth points.  Let $M_\h$ be a closed  hyperbolic $3$-manifold.  If $X$ and $M_\h$ are homeomorphic then
$\text{Vol} (X) \ge \text{Vol} (M_\h).$
\normalfont
\vskip 6pt

{\noindent}These inequalities do not follow from the Schl{\"a}fli formula for polyhedra \cite[pg.71]{CHK}.  Applying the Schl{\"a}fli formula requires a one-parameter family of cone-manifolds connecting $X$ to $M_\h$.  Such a family does not exist in general (see the end of Section \ref{cone manifolds}).

\vskip 6pt
The author would like to thank his advisor, Richard Canary, for his absolutely essential assistance.  From asking the initial question to editing the incorrect drafts, his advice was crucial at every stage of this research.  The author thanks Ian Agol for pointing out the application of Theorem \ref{technical theorem} to cone-manifolds.  Finally, the author thanks the referee for his excellent comments.

\subsection{Sketch of proof}
The powerful tool used here is a computation of spherical volume by Besson-Courtois-Gallot.  Let $Y$ be a compact manifold with universal cover $\widetilde{Y}$.  Philosophically, the spherical volume of $Y$ is the minimal $n$-dimensional Hausdorff measure of ``embeddings'' $Y \longrightarrow L^2 (\widetilde{Y})$.  (Obviously this is not the actual definition, but it captures the idea.)  In \cite{BCGlong}, Besson-Courtois-Gallot proved the spherical volume of a closed hyperbolic manifold is a dimensional constant times its volume.  Let $X$ be a compact Alexandrov metric space with curvature bounded below by $-1$.  (e.g. $X$ is a cone-manifold or $DC_M$.)  Generalizing a theorem of \cite{BCGlong}, we prove the spherical volume of $X$ is not greater than a dimensional constant times Vol$(X)$.  We then prove spherical volume is a bilipschitz invariant.  If $X$ is bilipschitz to a closed hyperbolic manifold $M_\h$, then these fact prove the inequality of Theorem \ref{technical theorem}, namely that Vol$(X) \ge \text{Vol}(M_\h)$. 

To use this inequality to study convex cores $C_M$, we must first change $C_M$ into a closed manifold.  This is accomplished by simply doubling $C_M$ across its boundary to obtain closed topological manifold $DC_M$.  (For technical reasons, a neighborhood of $C_M$ is doubled, but this is a detail.)  Second, $C_M$ is proven to be an Alexandrov space with curvature bounded below by $-1$. Next, we use the assumption that there exists $M_0 \in \text{cc}_{\text{top}}(N)$ such that $\partial C_{M_0}$ has totally geodesic boundary.  This assumption is used to find a bilipschitz map $DC_M \longrightarrow DC_{M_0}$.  $DC_{M_0}$ is a closed hyperbolic manifold.  Therefore, Theorem \ref{technical theorem} can be applied, yielding Theorem \ref{main2}.

Obtaining the inequality for cone-manifolds is simpler.  It was proven in \cite{BGP} that cone-manifolds with cone angles $\le 2\pi$ are Alexandrov spaces.  In Theorem \ref{cone1}, we assume a bilipschitz map to a closed hyperbolic manifold exists.  The desired volume inequality then follows immediately from Theorem \ref{technical theorem}.  For $3$-dimensional cone-manifolds, hard classical topology can be employed to promote a homeomorphism to a bilipschitz map basically for free, yielding Corollary \ref{cone2}.

\section{Preliminaries}
\label{Preliminaries}
The following is a review of the necessary definitions.

\subsection{Convex Core}
Let $M$ be a complete hyperbolic manifold.  Let $S \subseteq M$ be the union of all closed geodesics in $M$.  The \textit{convex core}, $C_M$, is the smallest closed convex subset of $M$ which contains $S$, in other words it is the closed convex hull of $S$ in $M$.  The convex core may also be defined as the smallest closed convex subset of $M$ such that the inclusion map is a homotopy equivalence.  

For finite volume hyperbolic manifolds, the convex core is the entire manifold.  Thus this is a useful object only in the infinite volume case, where $C_M$ is the smallest submanifold which carries all the geometry of $M$.

\subsection{Convex Cocompact}
A complete hyperbolic manifold $M$ is \textit{convex cocompact} if $C_M$ is compact.  These are the best behaved infinite volume hyperbolic spaces.  There is a natural deformation space associated with convex cocompact manifolds.  Fix a compact smooth manifold $N$ (usually with boundary).  Define $\text{cc}_{\text{top}} (N)$ to be the set of isometry classes of complete convex cocompact hyperbolic manifolds $M$ diffeomorphic to the interior of $N$.

For $3$-manifolds, the work of Thurston yields precise topological conditions on $N$ which imply $\text{cc}_{\text{top}} (N)$ is nonempty.  Specifically, $\text{cc}_{\text{top}} (N)$ is nonempty if and only if $N$ is a compact irreducible atoroidal $3$-manifold such that $\partial N$ is a nonempty collection of surfaces with negative Euler characteristic \cite{Mo}.

\subsection{Acylindrical}
Define $A := S^1 \times [0,1]$ to be a closed annulus.  A boundary preserving map $f: (A, \partial A) \longrightarrow (N, \partial N)$ is \textit{essential} if $f$ is $\pi_1$-injective and $f$ is not homotopic rel boundary to a map $g: A \longrightarrow \partial N$.  $N$ is \textit{acylindrical} if $\partial N$ is incompressible and there does not exist an essential map $f:(A, \partial A) \longrightarrow (N, \partial N)$.

\subsection{Alexandrov Spaces}
Let $Y$ be a complete locally compact geodesic metric space of finite Hausdorff dimension.  Consider points $p,\ q,\ r \in Y$.  By assumption there exist (not necessarily unique) geodesic paths connecting any pair of these points.  Any geodesic segment between $p$ and $q$ will be notated simply by $pq$ and its length by $|pq|$.  $\triangle pqr$ will denote a geodesic triangle formed by geodesic segments $pq,\ pr,$ and $ rq$.  $\widetilde{\triangle} pqr$ will denote the comparison triangle in $\mathbb{H}^2$ with side lengths $|pq|,\ |pr|$, and $ |rq|$.  Its vertices will be labelled in the obvious way by $\tilde{p} , \ \tilde{q}$, and $ \tilde{r}$. 

$Y$ is an \emph{Alexandrov space with curvature bounded below by $-1$} if (in addition to the above conditions) for some neighborhood $U_y$ of each point $y \in Y$ the following condition is satisfied:  For any triangle $\triangle{pqr}$ with vertices in $U_y$ and any point $s$ on the side $qr$ the inequality $|ps| \ge |\tilde{p} \tilde{s} |$ is satisfied, where $\tilde{s}$ is the point on the side $\tilde{q} \tilde{r}$ of the triangle $\widetilde{\triangle} pqr$ corresponding to $s$, that is, such that $|qs| = |\tilde{q} \tilde{s}|$, $|rs| = | \tilde{r} \tilde{s} |.$

\subsection{Volume growth entropy}
Volume growth entropy is a fundamental geometric invariant which plays a crucial role in the techniques developed by Besson-Courtois-Gallot.  Let $X$ be a metric space of Hausdorff dimension $n$, $\widetilde{X}$ be the universal cover of $X$, and $\mathcal{H}^n$ be $n$-dimensional Hausdorff measure.  The \textit{volume growth entropy} of $X$ is the number
$$h(\widetilde{X}) :=  \limsup_{R \rightarrow \infty} \frac{1}{R} \log \mathcal{H}^n ( B_{\widetilde{X}} (x,R)),$$
where $x$ is any point in $\widetilde{X}$, and the ball $B_{\widetilde{X}} (x,R)$ is in $\widetilde{X}$.  The volume growth entropy as defined is independent of the choice of $x \in \widetilde{X}$.  The following theorem of Burago, Gromov, and Perelman will be vitally important for this paper.

\begin{thm}
\label{Perel'man}
\cite[pg.40]{BGP}
If $X$ is an Alexandrov space with curvature bounded below by $-1$ and Hausdorff dimension $n \in \mathbb{N}$, then the volume growth entropy of $X$ is less than or equal to the volume growth entropy of $\mathbb{H}^n$, namely
$$h(\widetilde{X}) \le h(\mathbb{H}^n) = n-1.$$
\end{thm}

\subsection{Generalized differentiable and Riemannian structures}
See \cite{OS}.  Let $X$ be a toplogical space, $\Omega \subseteq X$, $n \in \mathbb{N}$, and $0 \le r < 2.$  A family $ \{ (U_\phi, \phi )\}_{\phi \in \Phi}$ is called a \emph{$C^r$-atlas on $\Omega \subseteq X$} if the following hold:

(1)  For each $\phi \in \Phi$, $U_\phi$ is an open subset of $X$.

(2)  Each $\phi \in \Phi$ is a homeomorphism from $U_\phi$ into an open subset of $\mathbb{R}^n$.

(3)  $\{ U_\phi \}_{\phi \in \Phi}$ is a covering of $\Omega$.

(4)  If two maps $\phi, \psi \in \Phi$ satisfy $U_\phi \bigcap U_\psi \not= \emptyset$, then
$$\psi \circ \phi^{-1} : \phi(U_\phi \bigcap U_\psi) \longrightarrow 
        \psi (U_\phi \bigcap U_\psi)$$
is $C^r$ on $\phi(U_\phi \bigcap U_\psi \bigcap \Omega)$.
\\

A family $\{g_\phi \}_{\phi \in \Phi}$ is called a \emph{$C^{r-1}$-Riemannian metric} associated with a $C^r$-atlas $\{ (U_\phi, \phi) \}_{\phi \in \Phi}$ on $\Omega \subseteq X$ if the following hold:

(1)  For each $\phi \in \Phi$, $g_\phi$ is a map from $U_\phi$ to the set of positive symmetric matrices.

(2)  For each $\phi \in \Phi$, $g_\phi \circ \phi^{-1}$ is $C^{r-1}$ on $\phi(U_\phi \bigcap \Omega)$.

(3)  For any $x \in U_\phi \bigcap U_\psi, \phi,\psi \in \Phi$, we have
$$g_\psi (x) = [ d(\phi \circ \psi^{-1})(\psi(x))]^t g_\phi (x)
        [d (\phi \circ \psi^{-1}) (\psi(x))].$$

The entire reason for introducing this terminology is the following theorem.

\begin{thm}
\label{Otsu}
\cite{OS}
Let $X$ be an $n$-dimensional Alexandrov space.  Then there exists a subset $S \subset X$ of Hausdorff dimension $\le n-1$ (a set of singular points), a $C^1$-atlas on $X \setminus S$, and a $C^0$-Riemannian metric on $X \setminus S$ associated with the atlas such that

(1)  The maps $\phi : U_\phi \longrightarrow \mathbb{R}^n$ of the $C^1$-atlas are locally bilipschitz.

(2)  For any $x,y \in X \setminus S$ and $\varepsilon > 0$, $x$ and $y$ can be joined by a path in $X \setminus S$ of length less than $d(x,y) + \varepsilon$.

(3)  The metric structure on $X \setminus S$ induced from the Riemannian structure coincides with the original metric of $X$.

(4)  In particular, the Riemannian metric induces a volume element $d\text{vol}_X$ on $X \setminus S$.  The measure on $X$ obtained by integrating this element equals $n$-dimensional Hausdorff measure on $X$ ($S$ has zero measure).
\end{thm}

\begin{rmk}
Statements (1),(2), and (4) above are not found in the beginning of \cite{OS}.  (1) can be found on page 651, (2) on page 654, and (4) on page 657.  In (4), volume elements are used instead of forms to avoid orientation issues.  For an explanation of this terminology, see \cite[pg.351]{S}.
\end{rmk}

{\noindent}We may therefore unambiguously define Vol$(X) := \mathcal{H}^n (X).$

The standard formulation of Rademacher's theorem states that locally Lipschitz maps $\mathbb{R}^n \longrightarrow \mathbb{R}^m$ are differentiable almost everywhere \cite[Thm.7.3]{Ma}.  Here we will use a generalized Rademacher's theorem stating that locally Lipschitz maps from $\mathbb{R}^n$ into a separable Hilbert space are differentiable almost everywhere.  This fact can be assembled from more general propositions in \cite{BL}.  Namely, use Corollary 5.12 on page 107, Proposition 6.41 on page 154 (note that Hilbert space is reflexive), and Proposition 4.3 on page 84.  This version of Rademacher's theorem remains true for Alexandrov spaces.

\begin{cor}
\label{Otsu corollary}
Let $X$ be an $n$-dimensional Alexandrov space, $L^2$ a separable Hilbert space, and $F: X \longrightarrow L^2$ a locally Lipschitz map.  Then $F$ is differentiable almost everywhere in $X$.
\end{cor}
\begin{pf}
By Theorem \ref{Otsu}, $F \circ \phi^{-1}$ is locally Lipschitz.  Therefore by Rademacher's theorem, it is differentiable a.e..  Points of differentiability in $\phi ( U_\phi \setminus S)$ are preserved under change of coordinates.
\end{pf}

\section{Spherical Volume}
\label{section spherical volume}
This section will introduce the notion of spherical volume, a geometric invariant defined by Besson-Courtois-Gallot~\cite{BCGold,BCGlong}.  It was used in the original proof of the Besson-Courtois-Gallot Theorem~\cite{BCGlong}, but later proofs do not make specific reference to it.  The computation of spherical volume stated in this section is nontrivial; it is the key fact used in this paper.

Let $X$ be a metric space of Hausdorff dimension $n$.  Then $L^2 (X)$ is the Hilbert space of square-integrable measurable functions on $X$ with respect to $n$-dimensional Hausdorff measure, and $S^\infty (X)$ denotes the unit sphere in $L^2 (X)$.  In a natural way, the isometry group of $X$ acts by isometries on $L^2 (X)$.  Namely for $\gamma \in \text{Isom}(X)$, $\gamma.f := f \circ \gamma^{-1}$.  This clearly restricts to an action on $S^{\infty} (X)$.

\begin{defn}
\label{vol}
Let $Y$ be an $n$-dimensional Alexandrov space with curvature bounded below.  (In particular, $Y$ could be a compact Riemannian manifold.)  Let $\widetilde{Y}$ be the universal cover of $Y$.  Let $\Theta: \widetilde{Y} \longrightarrow S^\infty (\widetilde{Y})$ be a $\pi_1 (Y)$-equivariant, positive, Lipschitz map.  For all points $x \in \widetilde{Y}$ where $\Theta$ is differentiable, define a ``metric'' $g_\Theta$ by
$$g_\Theta (u,v)_x := \langle d\Theta_x (u), d\Theta_x (v) \rangle_{L^2 (\widetilde{Y})},$$
where $u,v \in T_x \widetilde{Y}$.  As $\Theta$ is Lipschitz, Corollary \ref{Otsu corollary} implies $g_\Theta$ is defined almost everywhere.
\end{defn}

\begin{defn}
Using the previous notation, since $\Theta$ is assumed to be equivariant, $g_\Theta$ descends to a ``metric'' on $Y$ itself.  For $x \in Y$, if $T_x Y$ exists and has an orthonormal basis $\{ e_i \}$, then define 
$$\dt_Y (g_\Theta)(x) := \det (g_\Theta (e_i, e_j)_x)_{ij}.$$
$\dt_Y (g_\Theta)$ is defined a.e. and is, up to sign, independent of the choice of $\{e_i \}$.  Therefore it makes sense to define
$$Vol (\Theta) := \int_Y |\dt_Y (g_\Theta)|^{1/2} d\text{vol}_Y.$$
\end{defn}

\begin{defn} \cite{BCGold}
\label{SphereVol}
For $Y$ a compact Alexandrov space with curvature bounded below, define the set $\mathcal{L}$ to be all Lipschitz, positive, $\pi_1 (Y)$-equivariant maps from $\widetilde{Y}$ to the unit sphere $S^\infty (\widetilde{Y}) \subset L^2 (\widetilde{Y})$.  Define the \textit{spherical volume} of $Y$ to be
$$\text{SphereVol} (Y) := \inf_{\Theta \in \mathcal{L}} \{ Vol (\Theta) \}.$$
\end{defn}

\begin{thm}
\label{spherical volume}
\cite{BCGlong} If $M_\h$ is a closed oriented hyperbolic $n$-manifold, then
$$\text{SphereVol} (M_\h) = \left( \frac{(n-1)^2}{4n} \right)^{n/2}
         \text{Vol} (M_\h).$$
\end{thm}

\begin{rmk}
To the author's knowledge, there does not exist a published proof of Theorem \ref{spherical volume} exactly as it is stated here.  Theorem \ref{spherical volume} is proven in \cite{BCGlong}, but for a slightly different definition of spherical volume.
However, the theorem as stated can be easily assembled from published facts.  In \cite[pg.432]{BCGold}, using Definition \ref{SphereVol}, it is proven that 
$$\text{SphereVol} (M_\h) \le \left( \frac{(n-1)^2}{4n} \right)^{n/2}
         \text{Vol} (M_\h).$$
And in \cite[pg.744]{BCGlong}, the slightly modified definition of spherical volume is shown to be less than or equal to that given in Definition \ref{SphereVol}.  Together, this proves Theorem \ref{spherical volume}.  See also \cite[pg.627]{BCG3}.
\end{rmk}

\vskip 12pt
\section{The Analysis}
\label{main section}
The goal of this section is to prove Theorem \ref{technical theorem}, which will largely follow from Propositions \ref{propA} and \ref{propB}.  To begin, we extend an estimate proven in \cite[Prop.3.4]{BCGlong} to Alexandrov spaces.

\begin{prop}  \label{propA} 
Let $X$ be an Alexandrov space with curvature bounded below.  Then
$$\text{SphereVol}(X) \le \left( \frac{h(\widetilde{X})^2}{4n} \right)^{n/2} \text{Vol}(X).$$
\end{prop}

\begin{pf}By Theorem \ref{Otsu}, let $\Omega \subseteq \widetilde{X}$ be a set of full Hausdorff measure which has a $C^1$-differentiable structure and is a $C^0$-Riemannian manifold.  Let $L^2 (\widetilde{X})$ denote the Hilbert space of measurable real functions on the universal cover of $X$ with respect to $n$-dimensional Hausdorff measure.  Recall that $h(\widetilde{X})$ is the volume growth entropy of $(X,d)$.  For $c > h(\widetilde{X}) / 2$ define a map
$$\Psi_c: \widetilde{X} \longrightarrow L^2(\widetilde{X}) \ \ \text{by} \ \ [\Psi_c (x) ] (y) = e^{-c d(x,y)},$$ 
where d(x,y) denotes the lifted metric on $\widetilde{X}$.  It is an elementary estimate to show $c > h(\widetilde{X}) /2$ implies $\Psi_c(x) \in L^2 (\widetilde{X})$.  (Use that for $R \gg 0,\ \mathcal{H}^n (B(x,R)) \le e^{(h(\widetilde{X}) + \delta)R}$.)

The goal here is to use the map $\Psi_c$ to obtain the estimate of Proposition \ref{propA}.  To do so, we must first show $\Psi_c$ is Lipschitz, positive, and $\Gamma := \pi_1 (X)$ equivariant.  Positivity is obvious.  Recall that for $f \in L^2  (\widetilde{X})$ and $\gamma \in \Gamma$ the $\Gamma$-action on $L^2 (\widetilde{X})$ is given by
$$(\gamma . f)(x) = f ( \gamma^{-1} x) .$$
By this definition $\Psi_c$ is clearly $\Gamma$-equivariant.  Namely,
$$[ \Psi_c ( \gamma x) ] (y) = e^{-c d(\gamma x,y)}
    = e^{-c d(x, \gamma^{-1} y)} = [\gamma . \Psi_c (x)] (y).$$

\begin{lem}
\label{lipschitz}
$\Psi_c$ is Lipschitz.
\end{lem}
\begin{pf}
Pick points $x,y \in \widetilde{X}$.  The goal is to control
$$\int_{\widetilde{X}} | e^{-c d(x, \zeta)} - e^{-c d(y, \zeta)} |^2 d\zeta .$$

By the triangle inequality,
\begin{eqnarray*}
e^{-cd(x, \zeta)} - e^{-cd(y, \zeta)} & \le &
e^{-cd(x, \zeta)} - e^{-cd(y, x)} e^{-cd(x, \zeta)}  \\
    & = & e^{-cd(x, \zeta)} (1-e^{-cd(y,x)}).  \\
\Longrightarrow |e^{-cd(x,\zeta)}-e^{-cd(y, \zeta)}| & \le &
(e^{-cd(x,\zeta)} + e^{-cd(y, \zeta)} ) (1-e^{-cd(x,y)}).
\end{eqnarray*}

Therefore,
\begin{eqnarray*}  
\lefteqn{\int_{\widetilde{X}} | e^{-c d(x, \zeta)} - 
    e^{-c d(y, \zeta)} |^2 d\zeta }  \\
& \le & (1-e^{-cd(x,y)})^2  \int_{\widetilde{X}} [ e^{-c d(x, \zeta)} + e^{-c d(y, \zeta)} ]^2 d\zeta  \\
& \le & (1-e^{-cd(x,y)})^2 \int_{\widetilde{X}} [e^{-2c d(x, \zeta)} + 
    2 e^{-c d(x, \zeta)} e^{-c d(y, \zeta)} +
            e^{-2c d(y, \zeta)}] d\zeta  \\
& = & (1-e^{-cd(x,y)})^2 \left[ \| \Psi_c (x) \|^2 + \int_{\widetilde{X}}
    2 e^{-c d(x, \zeta)} e^{-c d(y, \zeta)} d\zeta + 
     \| \Psi_c (y) \|^2 \right] \\
& \le & (1-e^{-cd(x,y)})^2 \left[ \| \Psi_c (x) \|^2 +
     2 \| \Psi_c (x) \| \| \Psi_c (y) \| +  \| \Psi_c (y) \|^2 \right] . \\
\end{eqnarray*} 
Apply Cauchy-Schwarz to obtain the final line.

For all $\gamma \in \Gamma$
$$\| \Psi_c (\gamma x) \| = \| \gamma . \Psi_c (x) \| = \| \Psi_c (x) \|.$$
This means the real function
 $ \| \Psi_c \| : \widetilde{X} \longrightarrow \mathbb{R} $ descends to a
continuous function from compact $X$ to $\mathbb{R}$, forcing it to have a finite supremum $b < \infty$.  Therefore,
$$ \int_{\widetilde{X}} | e^{-c d(x, \zeta)} - e^{-c d(y, \zeta)} |^2 d\zeta
   \le (1-e^{-cd(x,y)})^2 (4b^2)$$
To finish things off, notice that for $t \ge 0$, $(1-e^{-ct}) \le ct$.  Applying this yields the final inequality
$$ \int_{\widetilde{X}} | e^{-c d(x, \zeta)} - e^{-c d(y, \zeta)} |^2 d\zeta
   \le (c d(x,y))^2 (4b^2),$$
which implies that
$$ \frac{ \| \Psi_c (x) - \Psi_c (y) \|}{d(x,y)} \le 2bc.$$
\end{pf}

We must compute the derivative of $\Psi_c$.  It is easy to determine what its derivative should be.  Namely, for a.e. $(x,y)$ and $v \in T_x \Omega$, naively differentiating yields
$$[ (d\Psi_c)_x (v)](y) = -c e^{-c d(x,y)} Dd_{(x,y)}(v,0).$$
(Dd notates the derivative of the metric.  The metric is $1$-Lipschitz and therefore differentiable a.e..)  This formula can be justified by a straightforward application of the Lebesgue dominated convergence theorem.  (Simply write down the definition of the derivative.  It involves a limit outside of an integral.  Each term of the integrand is dominated by a constant times $e^{-2c d(x,y)}$, which is integrable.  Therefore the limit may be pushed inside the integral.)

Recall that $S^\infty (\widetilde{X})$ denotes the unit sphere in $L^2 (\widetilde{X})$.  Define the map
$\pi : L^2 (\widetilde{X}) \setminus \{ 0 \} \longrightarrow S^{\infty} (\widetilde{X})$ to be radial projection.  If $\Psi_c$ is differentiable at $x$ then clearly $\pi \circ \Psi_c$ is still differentiable at $x$.  For $f, \theta \in L^2 (\widetilde{X})$, a two line computation yields
$$d \pi_f (\theta) = \frac{\theta}{\| f \|} - \frac{f}{\| f \|^3} \langle f, \theta \rangle.$$
Therefore, for a.e. $x \in \Omega$ and $v \in T_x \widetilde{X}$ we have 
$$d(\pi \circ \Psi_c)_x (v) = \frac{(d\Psi_c)_x (v)}{\| \Psi_c (x)\|} -
    \frac{\Psi_c (x)}{\| \Psi_c (x) \|^3} \langle \Psi_c (x), (d\Psi_c)_x (v) \rangle.$$
So finally, (as in \cite[pg.430]{BCGold})
\begin{eqnarray*}
\| d(\pi \circ \Psi_c)_x (v) \|^2 &=& \frac{\| (d\Psi_c)_x (v) \|^2}{\| \Psi_c (x) \|^2} - \frac{1}{\| \Psi_c (x) \|^4} \langle \Psi_c (x), (d \Psi_c)_x (v) \rangle^2. \qquad \quad (1)
\end{eqnarray*}

For simplicity, define $\Phi_c := \pi \circ \Psi_c$.  For all points 
$x \in \Omega$ at which $\Phi_c$ is differentiable, define the ``metric'' $g_\Phi$ by
$$g_\Phi (u,v) = \langle d\Phi_c (u) , d\Phi_c (v) \rangle_{L^2},$$
where $u,v \in T_x \Omega$.  (For brevity, we write $g_\Phi$ instead of $g_{\Phi_c}$.)

If $g_\Phi$ is defined at $x$, then pick an orthonormal basis $\{ e_i \} \subset T_x \Omega$.  The above computations yield the key inequality
\begin{eqnarray*}
\text{Trace}_\Omega (g_\Phi ) & =&  \sum_{i=1}^n \| (d \Phi_c)_x (e_i) \|^2
     \le  \sum \frac{\| (d\Psi_c)_x (e_i) \|^2}{\| \Psi_c (x) \|^2} 
            \ \ \ \ \text{ (by (1))}\\
    & = & \frac{1}{\| \Psi_c (x) \|^2} \int \sum \left(
        [(d\Psi_c)_x (e_i)] (y) \right)^2 dy \\
    &=& \frac{1}{\| \Psi_c (x) \|^2} \int \sum \left(
    -c e^{-c d(x,y)} \ Dd_{(x,y)}(e_i,0) \right)^2 dy \\
    &= & \frac{c^2}{\| \Psi_c (x) \|^2} \int e^{-2c d(x,y)}
    \left( \sum (Dd_{(x,y)}(e_i,0))^2 \right) dy.
\end{eqnarray*}
The distance metric $d$ is $1$-Lipschitz.  Thus for fixed $y$, the function $x \mapsto d(x,y)$ has an almost everywhere defined gradient vector field with norm bounded by $1$.  In other words,
$$\sum_{i=1}^n (Dd_{(x,y)} (e_i,0))^2 \le 1 \quad \text{for a.e. } x.$$
(In fact equality holds a.e..  This will not be used here.)  Combining this with the above inequality yields
$$\text{Trace}_\Omega (g_\Phi) \le c^2.$$

$\Phi_c$ is $\Gamma$-equivariant.  $\Gamma$ acts by isometries.  Therefore $g_\Phi$ is invariant under the action of $\Gamma$ and $\dt_\Omega (g_\Phi)$ descends to a function on $X$ defined a.e..  Since geometric mean is dominated by arithmetic mean, we have the inequality
$$\sqrt{| \text{det}_\Omega ( g_\Phi) |} 
\le (\frac{1}{n} \text{Trace}_\Omega (g_\Phi))^{n/2} \le (c^2 / n)^{n/2}.$$
Since $c> h(\widetilde{X})/2$ is arbitrary, this completes the proof of Proposition \ref{propA}. \end{pf}

\vskip 12pt

Spherical volume is in fact a bilipschitz invariant.
\begin{prop}\label{propB}
Let $F: X \longrightarrow Y$ be a bilipschitz homeomorphism between $n$-dimensional Alexandrov spaces with curvature bounded below.  Then 
$$\text{SphereVol}(X) = \text{SphereVol}(Y).$$
\end{prop}
\begin{pf}
Let $F$ also denote the lifted bilipschitz map $\widetilde{X} \longrightarrow \widetilde{Y}$.  Let $\Phi: \widetilde{Y} \longrightarrow S^\infty (\widetilde{Y})$ be a Lipschitz, positive, $\Gamma$-equivariant map.  Consider the composition $\Phi \circ F: \widetilde{X} \longrightarrow S^\infty (\widetilde{Y}).$  Use $\Phi \circ F$ to define the ``metric'' $g_{\Phi \circ F}$ on $X$.  A computation shows that 
$$\dt_{X} (g_{\Phi \circ F}) (p) = (\text{Jac}F)(p)^2 \cdot
     \dt_{Y} (g_{\Phi})(F(p)) \text{  a.e..}$$
Therefore
\begin{eqnarray*}
\int_{X} | \dt_{X} (g_{\Phi \circ F}) |^{1/2} d\text{vol}_{X} &=& \int_{X} |\text{Jac}F|  \cdot | \dt_{Y} (g_{\Phi}) \circ F |^{1/2} d\text{vol}_{X} \\
&=& \int_Y |\dt_Y (g_{\Phi})|^{1/2} d\text{vol}_Y \le \text{SphereVol}(Y).
\end{eqnarray*}

The map $\mathcal{I} : L^2 (\widetilde{Y}) \longrightarrow L^2 (\widetilde{X})$ defined by 
$$\mathcal{I} : f \longmapsto (f \circ F) \cdot |\text{Jac}F |^{1/2}$$
is a $\Gamma$-equivariant isometry taking positive functions to positive functions.  By composing this isometry with $\Phi \circ F$ we obtain a Lipschitz, positive, $\Gamma$-equivariant map
$$\mathcal{I} \circ \Phi \circ F : \widetilde{X} \longrightarrow 
        S^\infty (\widetilde{X}).$$
Since $\mathcal{I}$ is an isometry, we obtain
$$\text{SphereVol}(X) \le Vol(\mathcal{I} \circ \Phi \circ F) = \int_{X} | \dt_{X} (g_{\Phi \circ F}) |^{1/2} d\text{vol}_{X} \le  \text{SphereVol}(Y).$$
The opposite inequality is proven by reversing the roles of $Y$ and $X$.
\end{pf}

We are now ready to prove
\begin{thm}
\label{technical theorem}
Let $X$ be an Alexandrov space with curvature bounded below by $-1$, $M_\h$ a closed  hyperbolic manifold.  If $X$ and $M_\h$ are bilipschitz, then
$$\text{Vol}(X) \ge \text{Vol}(M_\h).$$
\end{thm}

\begin{pf}
If $M_\h$ is nonorientable then lift both $X$ and $M_\h$ to an oriented metric double cover.  If Theorem \ref{technical theorem} is true for these double covers, then it follows also for $X$ and $M_\h$.  Therefore assume without a loss of generality that $M_\h$ is oriented.

Recall that by Theorem \ref{Perel'man}, $h(\widetilde{X}) \le (n-1)$.  By Theorem \ref{spherical volume}, 
$$\text{SphereVol}(M_\h) =  \left( \frac{(n-1)^2}{4n} \right)^{n/2} Vol(M_\h).$$
By combining these facts with Propositions \ref{propA} and \ref{propB} we obtain
\begin{eqnarray*}
\lefteqn{\left( \frac{(n-1)^2}{4n} \right)^{n/2} Vol(M_\h) = \text{SphereVol}(M_\h)} \\
  & &  = \text{SphereVol} (X) 
    \le  \left(\frac{h(\widetilde{X})^2}{4n} \right)^{n/2} \text{Vol}(X) 
    \le \left( \frac{(n-1)^2}{4n} \right)^{n/2} Vol(X).
\end{eqnarray*}
This completes the proof of Theorem \ref{technical theorem}.
\end{pf}

\section{Application to Infinite Volume Hyperbolic Manifolds}

Recall the definition of the invariant $\mathcal{V}$.

\begin{defn}
Let $N$ be a compact smooth  $n$-manifold.  Then
$$\mathcal{V} (N) := \inf_{M \in \text{cc}_{\text{top}}(N)} \{ \text{Vol} (C_M) \}.$$
\end{defn}

\begin{thm}
\label{main2}
Let $N$ be a smooth  $n$-manifold.  If there exists $M_0 \in \text{cc}_{\text{top}}(N)$ such that $\partial C_{M_0} \subset M_0$ is a totally geodesic submanifold, then 
$$\mathcal{V}(N) = \text{Vol} (C_{M_0}).$$
\end{thm}

Assume there exists such an $M_0 \in \text{cc}_{\text{top}}(N)$.  Pick $M \in \text{cc}_{\text{top}} (N)$.  To prove the theorem, it is enough to show
$$\text{Vol} (C_{M_0}) \le \text{Vol} (C_M).$$

As a first step towards applying Theorem \ref{technical theorem}, we will establish the required bilipschitz map.

\begin{lem}
\label{lipschitzeomorphism}
Let $\overline{\mathcal{N}_{\varepsilon} (C_M)}$ denote a closed $\varepsilon$-neighborhood of $C_M$.  Then $\overline{\mathcal{N}_{\varepsilon} (C_M)}$ and $C_{M_0}$ are bilipschitz.
\end{lem}
\begin{pf}
As the hyperbolic metrics in consideration are convex cocompact, $C_M$ and $C_{M_0}$ are both compact.  By virtue of its beautiful boundary, $C_{M_0}$ is diffeomorphic to $N$.  $\overline{\mathcal{N}_{\varepsilon} (C_M)}$ is a $\mathcal{C}^{1,1}$-manifold with boundary, $\mathcal{C}^{1,1}$-diffeomorphic to $N$.  (This argument follows from the fact that the complement of any $\varepsilon$-neighborhood of the convex core is a $\mathcal{C}^{1,1}$-smooth product of a surface and an interval \cite{EM}.)  Therefore $\overline{\mathcal{N}_{\varepsilon} (C_M)}$ and $C_{M_0}$ are $\mathcal{C}^{1,1}$-diffeomorphic.  As they are both compact, the diffeomorphism is bilipschitz.
\end{pf}

Define $X^\varepsilon$ and $X_0$ to be the metric doublings of $\overline{\mathcal{N}_{\varepsilon} (C_M)}$ and $C_{M_0}$ across their respective boundaries.  $X_0$ is a closed hyperbolic manifold.  Clearly the bilipschitz homeomorphism of the lemma can be doubled to a bilipschitz homeomorphism $F: X_0 \longrightarrow X^\varepsilon$.  Define the subset
$$\Omega := X^\varepsilon \setminus \partial \overline{\mathcal{N}_{\varepsilon} (C)}.$$
$\Omega \subseteq X^\varepsilon$ is an open subset of full $n$-dimensional Hausdorff measure which is a Riemannian manifold.

To apply Theorem \ref{technical theorem}, it remains to prove that $X^\varepsilon$ is an Alexandrov space with curvature bounded below by $-1$, and that $X^\varepsilon$ has Hausdorff dimension $n$.  To do so we will use the following theorem.

\begin{thm}
\label{main1}
\cite{BGP,P}
Let $C$ be a closed strictly convex $n$-dimensional ($n \ge 2$) submanifold of a complete $n$-dimensional hyperbolic manifold.  Assume the boundary of $C$ is at least $\mathcal{C}^1$ smooth.  Let $X$ be the metric space obtained by doubling $C$ across its boundary.  $X$ is an Alexandrov space with curvature bounded below by $-1$.
\end{thm}

\begin{rmk}
A more general version of this theorem is stated without proof in \cite[pg.54]{BGP}.  A proof can be found in the unpublished manuscript \cite[pg.28]{P}.  To the author's knowledge, a published proof does not exist.  For completeness, an elementary proof of Theorem \ref{main1} (avoiding the beautiful machinery of Perelman \textit{et al.}) has been included in this paper as an appendix.
\end{rmk}

\vskip 6pt
\begin{lem}
Let $M$ be a complete hyperbolic n-manifold with convex core $C_M$.  The metric doubling of a closed $\varepsilon$-neighborhood of $C_M$ across its boundary is an Alexandrov space with curvature bounded below by $-1$.  Further, it has Hausdorff dimension $n$.
\end{lem}
\begin{pf}
Two standard facts of hyperbolic geometry are that the boundary of an $\varepsilon$-neighborhood of $C_M$ is $\mathcal{C}^{1,1}$ smooth, and the closed $\varepsilon$-neighborhood of $C_M$ is strictly convex \cite{EM}.  Theorem \ref{main1} can now be applied to show $D \overline{\mathcal{N}_\varepsilon (C_M)}$ is an Alexandrov space with curvature bounded below by $-1$.  

Now to show $D \overline{\mathcal{N}_\varepsilon (C_M)}$ has Hausdorff dimension $n$.  $\partial \overline{\mathcal{N}_\varepsilon (C_M)}$ is $\mathcal{C}^{1,1}$ and of topological dimension $n-1$, implying it is has $n$-dimensional Hausdorff measure zero.  $( D \overline{\mathcal{N}_\varepsilon (C_M)}\, \setminus \, \partial \overline{\mathcal{N}_\varepsilon (C_M)})$ is a Riemannian $n$-manifold.  Therefore $D \overline{\mathcal{N}_\varepsilon (C_M)}$ has Hausdorff dimension $n$.
\end{pf}

Therefore, applying Theorem \ref{technical theorem} yields
$$\text{Vol} (X^\varepsilon) \ge \text{Vol} (X_0).$$
Obtaining the inequality of Theorem \ref{main2} is now trivial.  For all 
$\varepsilon > 0$,
$$Vol(\overline{\mathcal{N}_{\varepsilon} (C_M)}) = (1/2) Vol(X^\varepsilon) \ge 
    (1/2) Vol(X_0) = Vol(C_{M_0}).$$
Therefore,
$$Vol(C_M) \ge Vol(C_{M_0}).$$
This completes the proof of Theorem \ref{main2}.
\square 

\vskip 12pt
As an immediate corollary, we can now relate $\mathcal{V} (N)$ to the Gromov norm of $DN$, the topological doubling of $N$ across its boundary.  Let $\| [\  \cdot \  ] \|$ denote Gromov norm.  For a definition of this norm, and a scintillatingly beautiful proof of the following theorem, see \cite{Th}.

\vskip 6pt
\noindent\textbf{Theorem} (Gromov)\textbf{:}  If $X$ is a closed oriented hyperbolic $n$-manifold, and $v_n$ is the volume of a regular ideal simplex in $\mathbb{H}^n$, then
$$\| [X] \| = \text{Vol}(X) / v_n.$$

\vskip 6pt
\begin{cor}
\label{Gromov}
Let $N$ be a smooth oriented $n$-manifold.  If there exists $M_0 \in \text{cc}_{\text{top}} (N)$ such that $\partial C_{M_0} \subset M_0$ is a totally geodesic codimension one submanifold, then 
$$\| [DN] \| = \frac{2 \mathcal{V}(N)}{v_n}.$$
\end{cor}
\begin{pf}
Simply double $C_{M_0}$ across its boundary to obtain a closed hyperbolic manifold diffeomorphic to $DN$.  Apply Gromov's theorem.
\end{pf} 

Like the Gromov norm, $\mathcal{V}$ behaves well under finite covers.

\begin{cor}
Let $N$ be a smooth  $n$-manifold.  Assume there exists $M_0 \in \text{cc}_{\text{top}} (N)$ such that $\partial C_{M_0} \subset M_0$ is a totally geodesic submanifold.  If $\widetilde{N}$ is a covering space of $N$ of degree $k<\infty$, then $\mathcal{V}(\widetilde{N}) = k \mathcal{V}(N)$.
\end{cor}
\begin{pf}
There exists an $\widetilde{M_0} \in \text{cc}_{\text{top}} (\widetilde{N})$ such that $\partial C_{\widetilde{M_0}}$ has totally geodesic boundary and $\widetilde{M_0}$ is a geometric degree $k$ cover of $M_0$.
\end{pf}

\vskip 12pt
\textbf{Application to 3-Manifolds}
\vskip 6pt

For $n > 3$, the geometry of hyperbolic $n$-manifolds is not well understood.  Little is known about what types of manifolds might satisfy the hypotheses of Theorem \ref{main2}.  For this reason, the most interesting results are obtained when attention is restricted to $3$-manifolds.  In a certain sense, a ``generic'' topological $3$-manifold will be acylindrical.  (Roughly, if the boundary components are sufficiently inter-tangled, then essential cylinders should not exist and the boundary should be incompressible.)  With this in mind the following consequence of Mostow Rigidity and Thurston's Hyperbolization Theorem shows that the hypotheses of Theorem \ref{main2} do apply to a large class of $3$-manifolds.

\vskip 6pt
\noindent\textbf{Corollary of Rigidity and Hyperbolization:}\itshape
~\cite[pg.14]{Th2} Let $N$ be a compact irreducible atoroidal $3$-manifold such that $\partial N$ is a nonempty collection of surfaces with negative Euler characteristic.  $N$ is acylindrical if and only if there exists a unique $M_0 \in \text{cc}_{\text{top}}(N)$ such that $\partial C_{M_0}$ is totally geodesic. \normalfont

\begin{rmk}
All facts in this section have been stated without parabolics.  With some notational effort, more general statements can be made.
\end{rmk} 

This remarkable corollary immediately suggests the conjecture mentioned in Section \ref{intro}.  Let $\mathcal{I}(N)$ denote the set of isometry classes of hyperbolic $3$-manifolds homotopy equivalent to $N$.  Recall the volume function Vol$: M \in \mathcal{I}(N) \longmapsto \text{Vol}(C_M).$

\begin{conj}
\label{mainconj}
Retain the above notation.  If $N$ is acylindrical, then for all $M\in \mathcal{I}(N) \setminus\{ M_0 \}$,
$$ Vol (C_{M_0}) < Vol (C_M).$$
In other words, $M_0$ is the \textit{unique global minimum} of the function Vol.
\end{conj}

{\noindent}The first partial answer to this conjecture was provided by Bonahon.  With entirely different techniques he proved the following.

\begin{thm}
\cite{Bon} Let $M_0 \in \text{cc}_{\text{top}} (N)$ be such that $C_{M_0}$ is codimension $0$ and $\partial C_{M_0}$ is totally geodesic.  Then $M_0$ is a \textit{strict local minimum} of the function Vol.
\end{thm}

{\noindent}In this context, Theorem \ref{main2} is another partial answer to
 Conjecture \ref{mainconj}. 

\begin{thm}
\label{3-cor}
Let $N$ be an acylindrical compact irreducible atoroidal $3$-manifold such that $\partial N$ is a nonempty collection of surfaces with negative Euler characteristic.  Then there exists an $M_0 \in \text{cc}_{\text{top}}(N)$ such that $\partial C_{M_0}$ is totally geodesic, and $M_0$ is a \textit{global minimum} of the function Vol over $\mathcal{I}(N)$.
\end{thm}
\begin{pf}
Pick an $M \in \mathcal{I}(N)$.  It suffices to show that $\text{Vol}(C_M) \ge \text{Vol}(C_{M_0})$.  For geometrically infinite manifolds, this inequality is trivial.  So let of first assume that $M$ is convex cocompact.  Then because $N$ is acylindrical, $M$ is in fact homeomorphic to int$(N)$ \cite[Lem.X.23,pg.235]{J}.  Thus the desired inequality follows from Theorem \ref{main2}.

Now assume that $M$ is geometrically finite, but not convex cocompact.  By \cite{BB}, $M$ is the strong limit of a sequence of convex cocompact manifolds.  Vol is continuous under strong limits \cite{Ta}.  This proves the desired inequality.
\end{pf}

A skeptic could accuse this theorem of proving one object we do not understand is equal to another object we do not understand.  This is not true.  Quite a bit is known about hyperbolic $3$-manifolds with totally geodesic boundary.  Most importantly, Kojima~\cite{Ko2} proved they can always be geometrically decomposed into partially truncated hyperbolic polyhedra.  Paoluzzi and Zimmermann~\cite{PZ} constructed an infinite family of such manifolds with one boundary component.  In some cases, the volume Vol$(C_{M_0})$ can be computed explicitly (i.e. actual numbers!) using these decompositions into truncated polyhedra.  For a list of such volumes, see \cite{Ush}.

\section{Application to Cone Manifolds}
\label{cone manifolds}

Theorem \ref{technical theorem} may also be applied to cone-manifolds with all cone angles $\le 2\pi$.

\begin{defn}
\cite[pg.53]{CHK}
An \emph{$n$-dimensional cone-manifold} is a manifold, $X$, which can be triangulated so that the link of each simplex is piecewise linear homeomorphic to a standard sphere and $X$ is equipped with a complete path metric such that the restriction of the metric to each simplex is isometric to a geodesic simplex of constant curvature $K$.  The singular locus $\Sigma$ consists of the points with no neighborhood isometric to a ball in a Riemannian manifold.
\end{defn}

It follows that

{\noindent}$\bullet \ \  \Sigma$ is a union of totally geodesic closed simplices of dimension $n-2$.

{\noindent}$\bullet \ $   At each point of $\Sigma$ in an open $(n-2)$-simplex, there is a \emph{cone angle} which is the sum of the dihedral angles of the $n$-simplices containing the point.

(Notice that cone-manifolds whose singular loci have vertices are allowed.)

\begin{lem}
\cite[pg.7]{BGP}
If all cone angles of an $n$-dimensional cone-manifold $X$ are $\le 2\pi$, and $K \ge -1$, then $X$ is an Alexandrov space with curvature bounded below by $-1$.
\end{lem}

An $n$-dimensional cone-manifold clearly has Hausdorff dimension $n$.  Therefore we have the following theorem.

\begin{thm}
\label{cone1}
Let $X$ be an $n$-dimensional cone-manifold with all cone angles $\le 2\pi$ and $K \ge -1$.  Let $M_\h$ be a closed  hyperbolic $n$-manifold.  If $X$ and $M_\h$ are bilipschitz then
$$\text{Vol} (X) \ge \text{Vol} (M_\h).$$
\end{thm}

In the case when $n=3$, applying the Manifold Hauptvermutung \cite{R} yields a stronger corollary.  To prove it, a bit of classical terminology is required.  A \emph{combinatorial manifold} is a simplicial complex $K$ such that the link of each simplex is piecewise linear homeomorphic to a standard sphere.  (Notice that a cone-manifold is a combinatorial manifold.)  A \emph{$\mathcal{C}^1$-triangulation} of a smooth manifold $Z$ is a combinatorial manifold $K$ together with a homeomorphism $f: K \longrightarrow Z$ which is a piecewise $\mathcal{C}^1$-diffeomorphism on each simplex.  A result of Whitehead \cite[pg.822, Thm.7]{Wh} states that any closed smooth manifold $Z$ admits a $\mathcal{C}^1$-triangulation.
 
\begin{cor}
\label{cone2}
Let $X$ be a $3$-dimensional cone-manifold with all cone angles $\le 2\pi$ and $K \ge -1$.  Let $M_\h$ be a closed hyperbolic $3$-manifold.  If $X$ and $M_\h$ are homeomorphic then
$$\text{Vol} (X) \ge \text{Vol} (M_\h).$$
\end{cor}
\begin{pf}
It is sufficient to upgrade the homeomorphism to a bilipschitz homeomorphism.  Pick a $\mathcal{C}^1$-triangulation $f: K \longrightarrow M_\h$.  Let $K$ have the obvious piecewise Euclidean metric.  For each simplex $\sigma$ of $K$, $f|_{\sigma}$ is clearly bilipschitz.  By compactness, $f$ is globally bilipschitz.

$K$ and $X$ are homeomorphic.  By \cite{M}, this implies there exists a piecewise linear homeomorphism $X \longrightarrow K$.  A piecewise linear homeomorphism between closed combinatorial manifolds is bilipschitz.  Therefore, by postcomposing with $f$, $X$ and $M_\h$ are bilipschitz.  Now the previous theorem can be applied.
\end{pf}

Previously, results similar to these could be obtained by using the Schl{\"a}fli formula for polyhedra \cite[pg.71]{CHK}.  But to compare cone-manifolds $X$ and $M_\h$ using the Schl{\"a}fli formula it is necessary to have a one-parameter family of cone-manifolds connecting them.  Such a path in deformation space is \emph{not} necessary to apply Corollaries \ref{cone1} and \ref{cone2}.  We now sketch an example where such a path does not exist.  Let $\gamma$ be a simple closed geodesic in a closed hyperbolic manifold $M_\h$.  There exists a simple closed curve $\gamma' \subset M_\h$ such that $\gamma'$ is homotopic to $\gamma$, $M_\h \setminus \gamma$ is not homeomorphic to $M_\h \setminus \gamma'$, and $M_\h \setminus \gamma'$ admits a hyperbolic structure \cite{My}.  There exist cone-manifold deformations of this hyperbolic structure with singular locus $\gamma'$ and strictly positive cone angle \cite[pg.99]{CHK}.  If it were possible to increase the cone angle all the way to $2 \pi$, then by Mostow Rigidity  $M_\h \setminus \gamma'$ would be homeomorphic to $M_\h \setminus \gamma$.  Therefore such a family of deformations does not exist.

\section{Concluding Remarks}
\label{conclusion}
There are a couple obvious ways these results could be improved.

\noindent\textbf{1.}  Prove that hyperbolic metrics on closed manifolds uniquely minimize volume over Alexandrov metrics in the same bilipschitz class.  This is probably impossible using only spherical volume, as spherical volume was insufficient for proving the uniqueness statement of the Besson-Courtois-Gallot theorem. 

\noindent\textbf{2.}  Control the growth of the function Vol as a convex cocompact hyperbolic manifold $M$ moves away from the minimum $M_0$ in $\text{cc}_{\text{top}} (N)$.  It is a bit hidden the way the proofs were written here, but this could be done by improving the estimate of $h(\widetilde{X})$.  In $3$ dimensions, it seems plausible that as the bending lamination on $C_M$ becomes more extreme, $h(\widetilde{X})$ should go down.  Proving this strong result would require knowing something about Alexandrov spaces: how the local geometry near singularities affects the large scale geometry.

\appendix
\section{Proof of Theorem \ref{main1}}
Recall the statement of the theorem.

\vskip 6pt
\noindent\textbf{Theorem \ref{main1}.}  \itshape Let $C$ be a closed strictly convex $n$-dimensional ($n \ge 2$) submanifold of a complete $n$-dimensional hyperbolic manifold.  Assume the boundary of $C$ is at least $\mathcal{C}^1$ smooth.  Let $X$ be the metric space obtained by doubling $C$ across its boundary.  $X$ is an Alexandrov space with curvature bounded below by $-1$. \normalfont
\vskip 6pt

This theorem will be proven using a slightly different definition of Alexandrov space.  Namely, $X$ will be proven to be an angled Alexandrov space with curvature bounded below by $-1$.  (For locally compact spaces, these notions are equivalent, see Remark \ref{equivalent}.)  By doing so, the arguments remain a bit closer to those used on more familiar geometric objects.  It is therefore necessary to define what will be meant by angle.

Let $Y$ be geodesic metric space.  $\widetilde{\angle} rpq$ will denote the angle in $\mathbb{H}^2$ between the sides $\tilde{p} \tilde{r}$ and $ \tilde{p} \tilde{q}$ of a comparison triangle $\widetilde{\triangle} pqr$.  Defining angles in $Y$ itself is a bit trickier.

\begin{defn}
Let $rpq$ denote the union of a geodesic segment $pr$ and a geodesic segment $pq$.  Let $\{ r_i \} \subset pr$ and $\{ q_i \} \subset pq$ be sequences of points not equal to $p$ such that 
$r_i \longrightarrow p$ and $q_i \longrightarrow p$.  Then define
\[\angle rpq := \lim_{i \rightarrow \infty} 
     \widetilde{\angle} r_i p q_i.  \]
Clearly one must show this limit exists in order for this definition to make any sense.
\end{defn}

\begin{defn}
An \textit{angled Alexandrov space with curvature bounded below by} $-1$ is a complete geodesic metric space $Y$ such that for all points $x \in Y$ there exists an open neighborhood $U_x$ of $x$ satisfying the following three conditions for all geodesic triangles $\triangle pqr$ with vertices in $U_x$:

(A)  The angle $\angle rpq$ is defined.

(B)  None of the angles of the triangle $\triangle pqr$ is less than the corresponding angle of the comparison triangle in $\mathbb{H}^2$, i.e. $\angle rpq \ge \widetilde{\angle} rpq$.

(C)  The sum of adjacent angles is equal to $\pi$, i.e. if $s$ is an interior point of a geodesic $pq$ then for any geodesic $sr$ we have $\angle psr + \angle rsq = \pi$.
\end{defn}

These conditions have simple geometric consequences for the space.  Condition (B) states that triangles in $Y$ must be fat, and condition (C) guarantees that geodesics do not branch.  

\begin{rmk}
\label{equivalent}
For $Y$ locally compact, $Y$ is an angled Alexandrov space with curvature bounded below by $-1$ if and only if it is an Alexandrov space with curvature bounded below by $-1$.  See \cite[pg.7]{BGP}, or \cite[pg.114]{BBI}.
\end{rmk}

Now to prove Theorem \ref{main1}.  Let $X$ be as in the statement of the theorem.  It is easy to see that $X$ is geodesic, complete, and locally compact.  It is therefore sufficient to prove $X$ satisfies conditions (A), (B), and (C).  Since these are all local conditions, nothing is lost by lifting to the universal cover of the ambient hyperbolic manifold in which $C$ lives.  Therefore assume without a loss of generality that $C$ is a simply connected subset of $\mathbb{H}^n$.

Let $C_1$ and $C_2$ denote the two isometric copies of $C$ imbedded into $X$.  As their boundaries coincide, denote their common boundary by $\partial C$.  Assume $C_2$ is the half with the opposite orientation of $C$.  For points $x \in X \setminus \partial C$ there exists a neighborhood $U_x \ni x$ entirely contained in either $int(C_1)$ or $int(C_2)$.  So we can assume $U_x$ is an open set of $\mathbb{H}^n$ and therefore trivially satisfies conditions (A), (B), and (C).  So the definition needs to be verified only for points lying on $\partial C$.

For points $p,q \in int(C_i)$ lying on the same side of $X$, convexity and negative curvature imply the geodesic $pq$ is unique and lies entirely in $C_i$.  This is true even if $p \in int(C_i)$ and $q \in \partial C$.  If $p,q \in \partial C$ then there exist exactly two geodesics connecting them.  One geodesic lies in $C_1$, the other lies in $C_2$, and they are exchanged by the natural isometric involution of $X$.  If $p \in int(C_1)$ and $q \in int(C_2)$ then by convexity of $C$ and minimality of geodesics any geodesic $pq$ (no longer necessarily unique) intersects the boundary $\partial C$ in a unique point $c$.  As $\partial C$ was assumed to be $\mathcal{C}^1$, it makes sense to speak of the tangent space to $\partial C$ at $c$.  Thus, the angles between the tangent space and the smooth geodesics $pc$, $qc$ are well defined.  The first lemma is to prove those angles of incidence are equal.  Informally speaking, we will show that for geodesics in $X$ intersecting $\partial C$ the angle of incidence equals the angle of reflection.  First this terminology must be clarified.

Consider the following model for the metric space $X$.  Embed $C_1$ and $C_2$ \textit{on top of each other} into $\mathbb{H}^n$ respectively by an orientation-preserving and an orientation-reversing isometry with identical range.  Let $\Phi$ be the gluing-together of these two isometries.  Let $c'$ be $\Phi (c)$, $p'$ be $\Phi (p)$, and $q'$ be $\Phi (q)$.  Then clearly the image of a geodesic segment $pq$ under $\Phi$ is a shortest path connecting $p'$ to $q'$ \textit{passing through} $\Phi( \partial C)$.  Let this shortest path be $\sigma$.  The analogy to keep in mind is that $\sigma$ is the path light would travel bouncing off of a reflective surface.

\begin{lem}
\label{lemma1}
\textit{(angle of incidence equals angle of reflection)}
Let $\vec{N}$ be the inward-pointing normal vector to $\Phi( \partial C)$ at $c'$.  Let $\angle_{\mathbb{H}^n}$ denote angles measured in $\mathbb{H}^n$.  Then
\[\angle_{\mathbb{H}^n} (\vec{N}, c'p')) =
         \angle_{\mathbb{H}^n} (\vec{N}, c'q') \]
More importantly, for any vector $\vec{v}$ in the tangent space of $C$ at $c'$ such that the inner product of $\vec{v}$ with $\vec{N}$ is nonnegative, the following inequality holds:
\[\angle_{\mathbb{H}^n} (\vec{N}, c'p')) +
         \angle_{\mathbb{H}^n} (\vec{N}, c'q')  \le
   \angle_{\mathbb{H}^n} (\vec{v}, c'p')) +
         \angle_{\mathbb{H}^n} (\vec{v}, c'q')   \le  \pi  \]
with equality on the right if and only if $\vec{v} \bot \vec{N}$.

\end{lem}

\begin{pf}
Define $\Pi$ to be the totally geodesic $(n-1)$-dimensional subspace of $\mathbb{H}^n$ tangent to $\Phi(X)$ at $c'$.  Let $\gamma$ be the shortest path in $\mathbb{H}^n$ (not necessarily contained in $\Phi (X)$) joining $p'$ to $q'$ passing through $\Pi$.  $\Pi$ and $\gamma$ intersect in a unique point.  If $\Pi \bigcap \gamma = c'$, then the hypotheses are trivially true by elementary geometry.  So assume $\Pi \bigcap \gamma$ is not $c'$.  See Figure 1.

\begin{figure}[ht]
\centering
\includegraphics{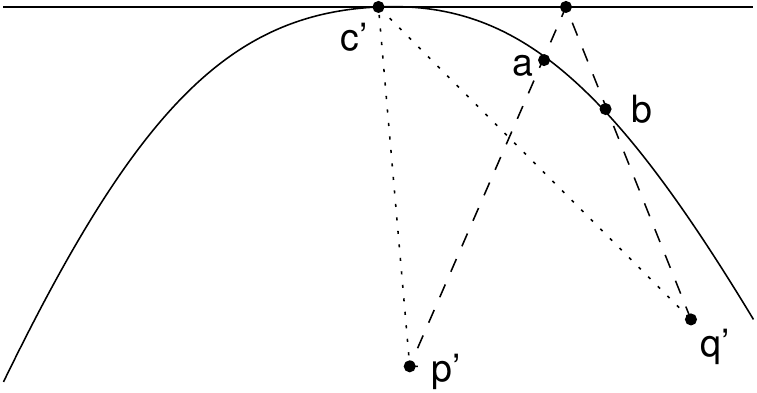}
\caption{the image of $X$ under $\Phi$}
\end{figure} 

By convexity of $C$, this implies that $\gamma$ intersects $\Phi(\partial C)$ in two points $a$ and $b$.  Recall $\sigma$ is the shortest path in $\Phi(X)$ connecting $p'$ to $q'$ passing through $\Phi(\partial C)$.  So by the triangle inequality
\[  \text{length of } \sigma \le d_{\mathbb{H}^n} (p' , a) + d_{\mathbb{H}^n} (a,b) + 
       d_{\mathbb{H}^n} (b, q')  \le  \text{length of} \ \gamma  . \]
Therefore the length of $\sigma$ is less than or equal to the length of $\gamma$.  But $\gamma$ was assumed to be minimal in $\mathbb{H}^n$.  This is a contradiction.  Therefore it follows that $\Pi \bigcap \gamma = c'$.

\end{pf}

This lemma provides a good picture for the behavior of geodesics in $X$ between $C_1$ and $C_2$.  What remains is to prove that angles are defined in $X$ and to find a simple formula for them.  The angle between any two geodesics leaving a point in the interior of either $C_1$ or $C_2$ is trivially defined and is equal to the corresponding angle in the image of $\Phi$ in $\mathbb{H}^3$.  Even for points on the boundary $\partial C$ the angle between two geodesics is obvious for geodesics both heading into the same side of $X$.  So the only interesting case is when two geodesic rays originating at $p \in \partial C$ head into opposite halves of $X$.  This means there exist $r \in int(C_1)$ and $q \in int(C_2)$ such that we are considering geodesics $pr$ and $pq$.  

Understanding this situation requires a different model.  Let 
\[\Psi_1: C_1 \longrightarrow 
\mathbb{H}^n \ \  \text{and} \  \ \Psi_2 :  C_2 \longrightarrow \mathbb{H}^n \] be \textit{orientation-preserving} isometries such that the images are tangent at their unique point of intersection
$$p' := \Psi_1 (p) = \Psi_2 (p).$$  
Let $S_1$ and $S_2$ be the half-spaces of $\mathbb{H}^n$ which contain $\Psi_1 (C_1)$ and $\Psi_2 (C_2)$ respectively and intersect in a plane through $p'$.  (See Figure 2.)

\begin{lem}
\label{lemma2}
The angle in $X$ formed by the geodesics $pr$ and $pq$ exists and is equal to the angle in $\mathbb{H}^n$ made by $\Psi_1 (pr)$ and $\Psi_2 (pq)$, namely
$$\angle qpr = \angle_{\mathbb{H}^n} \Psi_2(q) p' \Psi_1 (r).$$
\end{lem}
This lemma implies $X$ satisfies condition (A).

\begin{figure}[ht]
\centering
\includegraphics{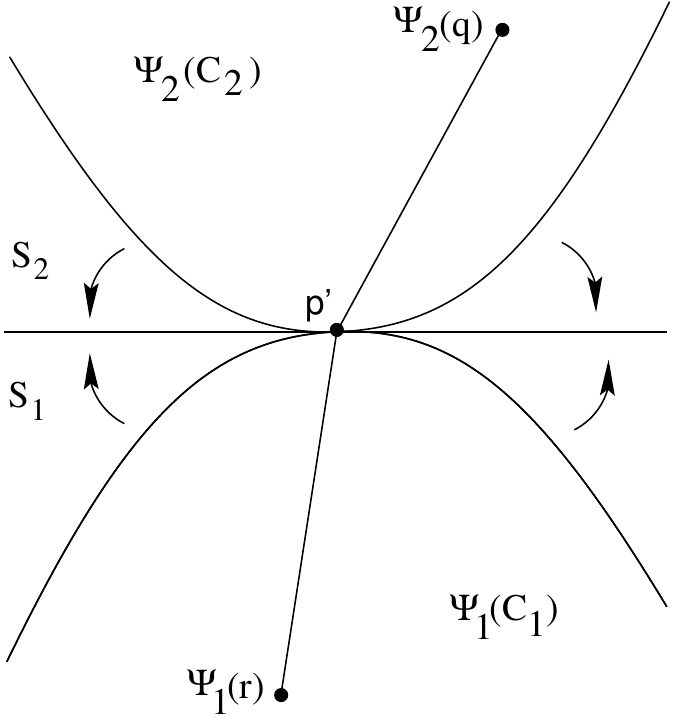}
\caption{the image of $X$ under $\Psi_1$ and $\Psi_2$}
\end{figure}  

\begin{pf}
By convexity of $C_i$ and the fact that $\Psi_1 (C_1)$ and $\Psi_2 (C_2)$ are tangent at $p'$, there exists a neighborhood $\mathcal{O} \subset \mathbb{H}^n$ of $p'$ and bilipschitz embeddings $f_i :\Psi_i (C_i) \bigcap \mathcal{O} \longrightarrow S_i$ such that:

(1)  There exists a function $K(R)$ such that $f_i$ is $K(R)$-bilipschitz on $B(p,R) \bigcap \Psi_i (C_i)$, and $\lim_{R \rightarrow 0} K(R) = 1$.

(2)  $f_1$ is the identity on geodesic ray $\Psi_1 (pr) \bigcap \mathcal{O}$.  $f_2$ is the identity on the geodesic ray $\Psi_2 (pq) \bigcap \mathcal{O}$.

(3)  The $f_i$ agree on $\Psi_i (\partial C_i)$ in such a way that they glue together to form a map $F$ from a neighborhood of $p$ in $X$ to $\mathbb{H}^3$.

The $f_i$ are obtained by spreading out $\Psi_i (C_i) \bigcap \mathcal{O}$ to fill up $S_i$ while preserving radial distances from $p'$.  Notice that $F$ also satisfies (trivial modifications of) properties (1) and (2).

Pick sequences $\{ r_i \} \subset int(pr)$ and $\{ q_i \} \subset int(pq)$ both converging to $p$.  By properties (1) and (2), the three-point metric spaces $\{ r_i, p, q_i \} \subset X$ and $\{ F(r_i) = \Psi_1 (r_i), p', F(q_i) = \Psi_2 (q_i) \} \subset \mathbb{H}^n$ are $K_i$-bilipschitz, with $K_i$ going to $1$.  Therefore
$$| \widetilde{\angle} q_i p r_i - \angle_{\mathbb{H}^n} F(q_i) p' F(r_i)| \longrightarrow 0.$$
We also know that
$$ \angle_{\mathbb{H}^n} F(q_i) p' F(r_i) = \angle_{\mathbb{H}^n} \Psi_2 (q_i) p' \Psi_1 (r_i)
    = \angle_{\mathbb{H}^n} \Psi_2 (q) p' \Psi_1 (r).$$
This implies that
$$ \widetilde{\angle} q_i p r_i \longrightarrow \angle_{\mathbb{H}^n} \Psi_2 (q) p' \Psi_1 (r).$$
This proves the desired equality.
\end{pf}

The easiest remaining condition to verify is condition (C).
\vskip 6pt
(C)  The sum of adjacent angles is equal to $\pi$, i.e. if $s$ is an interior point of a geodesic $pq$ then for any geodesic $sr$ we have $\angle psr + \angle rsq = \pi$.
\vskip 6pt
The only nontrivial case of (C) is for a geodesic segment straddling the boundary $\partial C$.

\begin{lem}
\label{conditionC}
Let $qr$ be a geodesic segment such that $r \in C_1$, $q \in C_2$, and $p$ is the unique point in the intersection $qr \, \cap \partial C$.  Let $s$ be a point of $X$.  This arrangement satisfies condition (C), i.e. $\angle rps + \angle spq = \pi$.
\end{lem}
\begin{pf}
Without a loss of generality, we may assume $s \in C_1$.  Consider the model used in Lemma \ref{lemma2}.  Specifically, map $X$ into $\mathbb{H}^n$ by 
$\Psi_1$ and $\Psi_2$ so that $\Psi_1(  C_1 )$ and 
$\Psi_2 (  C_2 )$ are tangent at the image of $p$.  Then by Lemma \ref{lemma1}, the geodesic segment $qr$ is mapped to a \textit{geodesic segment} in $\mathbb{H}^n$.  Therefore,
$$\angle_{\mathbb{H}^n} \Psi_2(q) \Psi_1(p) \Psi_1(s)  +
    \angle_{\mathbb{H}^n} \Psi_1(s) \Psi_1(p) \Psi_1(r)  = \pi .$$
Lemma \ref{lemma2} states that
$$ \angle_{\mathbb{H}^n} \Psi_2(q) \Psi_1(p) \Psi_1(s) = \angle qps.$$
It is obvious that
$$ \angle_{\mathbb{H}^n} \Psi_1(s) \Psi_1(p) \Psi_1(r) = \angle spr.$$
Therefore,
$$ \angle qps + \angle spr = \pi.$$
\end{pf}

What remains is to verify condition (B) for small triangles in $X$.  Recall
\vskip 6pt
(B)  None of the angles of the triangle $\triangle pqr$ is less than the corresponding angle of the comparison triangle in $\mathbb{H}^2$, i.e. $\angle rpq \ge \widetilde{\angle} rpq$.
\vskip 6pt
This requires the following formulation of a classical result of Alexandrov~\cite[pg.115]{BBI}. \\ \\
\textbf{Alexandrov's Lemma:  }\itshape Let points $a,b,c,d \in \mathbb{H}^2$ form a geodesic quadrilateral as in Figure 3.  Let $\Delta$ be the geodesic triangle in $\mathbb{H}^2$ with side lengths $|ab|, |bc|+|cd|, \text{ and } |ad|$.
\begin{figure}[ht]
\centering
\includegraphics{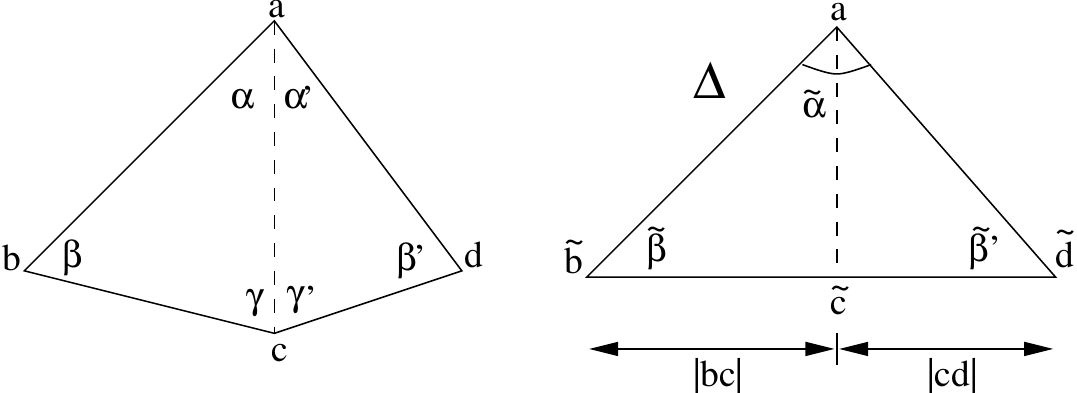}
\caption{}
\end{figure}  
If $\gamma + \gamma' \le \pi$, then
$\widetilde{\alpha} \ge \alpha + \alpha'$, $\widetilde{\beta} \le \beta$, $\widetilde{\beta'} \le \beta'$, and $d(\widetilde{a}, \widetilde{c}) \le d(a,c)$.
Moreover, if any of the above inequalities is an equality, then all the others are also equalities. \normalfont

\begin{lem}
\label{work1}
A triangle $\triangle qpr$ where $p \in \partial C, q \in int(C_1), $ and $r \in int(C_2)$ satisfies condition (B).
\end{lem}
\begin{pf}
Pick a geodesic $qr$ in $X$.  Let $c$ be the unique point where $qr$ intersects $\partial C$.  Consider the triangles $\triangle pqc \in C_1$ and $\triangle prc \in C_2$.  Let $Q$ be the (hyperbolic planar) quadrilateral formed by abstractly gluing $\triangle pqc$ to $\triangle prc$ along the edges (unfortunately labelled by identical notation) $pc \in C_1$ and $pc \in C_2$.  The idea is to compare $Q$ to the comparison triangle $\widetilde{\triangle} pqr$.  See Figure 4.  By Lemma \ref{lemma1}, we know $$\angle rcp + \angle qcp \le \pi.$$  As all side lengths are equal in $Q$ and $\widetilde{\triangle} pqr$, we can apply Alexandrov's lemma.  Therefore,
$$\angle prc \ge \widetilde{\angle} prc   \  \  \text{and} \ \ \angle pqc
 \ge \widetilde{\angle} pqc .$$

\begin{figure}[ht]
\centering
\includegraphics[scale=0.6]{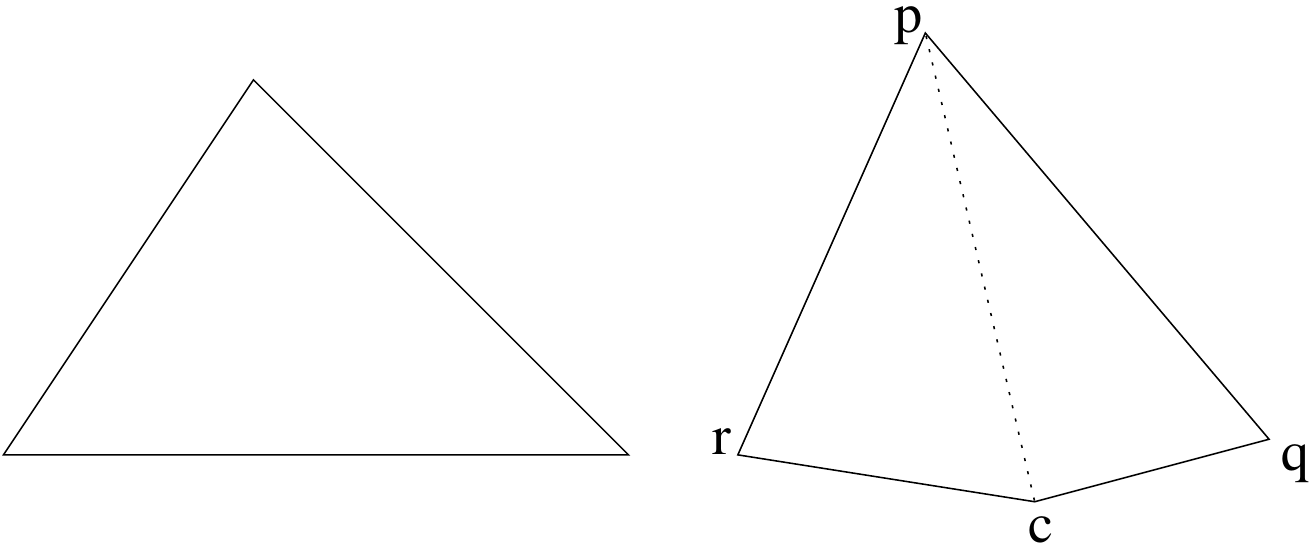}
\caption{the comparison triangle $\widetilde{\triangle} pqr$ and the glued-together quadrilateral $Q$}
\end{figure}  

All that remains is showing $\angle rpq \ge \widetilde{\angle} rpq .$  (Warning, Figure 4 is misleading in this case!)  This requires comparing the comparison triangle $\widetilde{\triangle} pqr$ to the triangle in $\mathbb{H}^n$ formed by $\Psi_1 (p), \ \Psi_1(q), \ \text{and} \ 
\Psi_2 (r)$ considered in Lemma \ref{lemma2}.  We know
\[ \angle rpq = \angle_{\mathbb{H}^n} \Psi_2(r) \Psi_1(p) \Psi_1(q) , \]
\[ d_{\mathbb{H}^n} (\Psi_1(p), \Psi_1(q)) = |pq| , \]
\[ d_{\mathbb{H}^n} (\Psi_1(p), \Psi_2(r)) = |pr| , \text{ and} \]
\[ d_{\mathbb{H}^n} (\Psi_1(q), \Psi_2(r)) \ge |qr| .\]
Therefore,
\[ \angle rpq =  \angle_{\mathbb{H}^n} \Psi_2(r) \Psi_1(p) \Psi_1(q)
    \ge \widetilde{\angle} rpq  .\]
\end{pf}

Now consider a triangle in $X$ which straddles the boundary $\partial C$ but has no vertices lying on $\partial C$.

\begin{lem}
A triangle $\triangle qpr$ where $p \in int(C_1), q \in int(C_1), $ and $r \in int(C_2)$ satisfies condition (B).
\end{lem}
\begin{pf}
This is a disappointingly messy case.  The idea is to cut $\triangle qpr$ into pieces for which the previous lemmas are valid.

Let $c$ be the unique point where $qr$ intersects the boundary $\partial C$. 
Let $Q$ be the quadrilateral formed by gluing the \textit{comparison triangle} $\widetilde{\triangle} rcp$ to $\triangle pcq$ along the appropriate edge.  (See Figure 5.)

\begin{figure}[ht]
\centering
\includegraphics[scale=0.6]{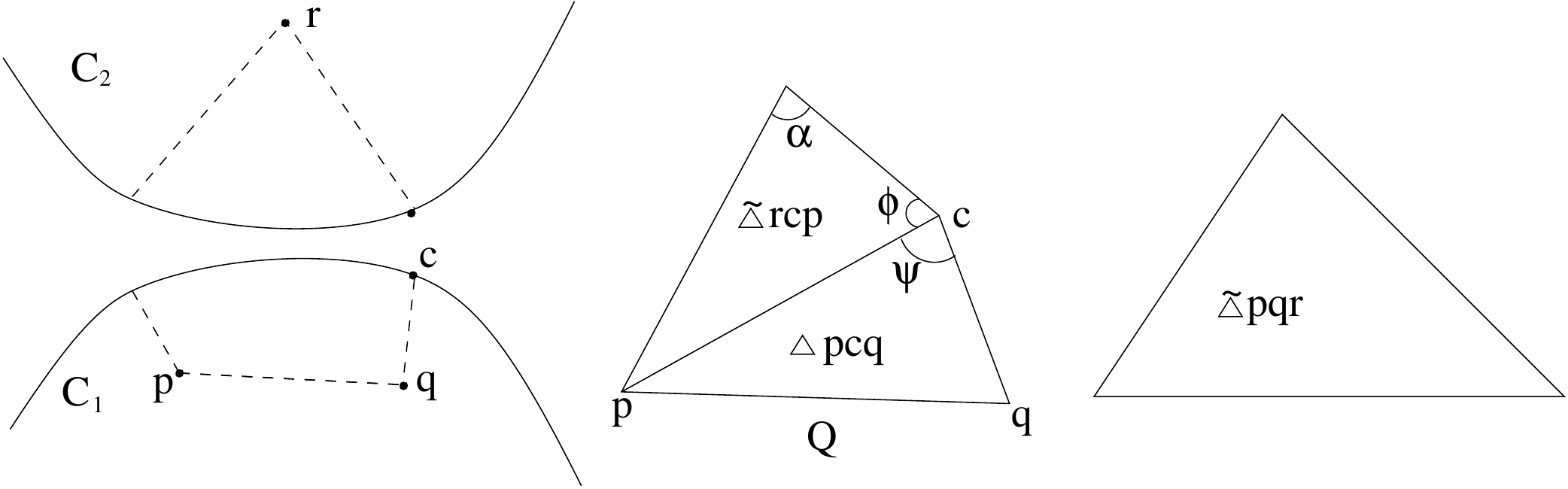}
\caption{a visual model for $\triangle qpr$, the quadrilateral $Q$, and comparison triangle $\widetilde{\triangle} pqr$ }
\end{figure}  

Since we could have just as easily cut along the other diagonal, by symmetry it is enough to show that
$$ \angle prq = \angle prc \ge \widetilde{\angle} prq, \text{ and} \ \ 
 \angle rqp = \angle cqp \ge \widetilde{\angle} rqp .  $$

Lemma \ref{work1} can be applied to $\triangle rcp$.  Therefore, 
Lemma \ref{work1} and condition (C) imply
$$\phi + \psi = \phi + \angle pcq \le \angle rcp + \angle pcq = \pi .$$

Using this we can apply Alexandrov's lemma to conclude $\widetilde{\angle} prq \le \alpha$.  By Lemma~\ref{work1}, $\alpha \le \angle prc = \angle prq$.  Therefore, $$\widetilde{\angle} prq \le \angle prq.$$  Finally, by again using Alexandrov's lemma we obtain
$$ \widetilde{\angle} rqp \le \angle cqp = \angle rqp .$$

\end{pf}

With this, we have proven that each $x \in \partial C$ has a neighborhood satisfying conditions (A), (B), and (C) where we let $U_x$ be simply $X$ in each case.  (This was possible only because we reduced to the case where $C$ is simply connected.)  This completes the proof of Theorem \ref{main1}.

\begin{sc}
\noindent
Department of Mathematics\\
University of Michigan\\
East Hall, 525 E University Ave\\
Ann Arbor, MI 48109\\
\end{sc}

\end{document}